\documentclass[10pt]{article}
\usepackage{amssymb,cite}
\usepackage{amsmath}
\usepackage{epsfig}
\usepackage{graphicx}
\graphicspath{ {./images/} }
\usepackage[all]{xy}
\usepackage{color}
\usepackage{enumitem}  
\usepackage{kotex}
\usepackage{authblk}
\usepackage{ulem} 
\usepackage{mathrsfs} 
\usepackage{float}

\usepackage[mathlines]{lineno} 
\usepackage{multirow}

\usepackage[color=green!40]{todonotes}
\usepackage{setspace}
\usepackage{amsthm,amsfonts,amssymb,epsfig,graphics,amsmath,amsbsy}
\usepackage{subfig}

\begin{document}
\def\mathbi#1{\textbf{\em #1}}
\newcommand{\intL}{\int\limits}
\newcommand{\half}{^\infty_0 }
\newcommand{\intR}{\int\limits_{\mathbb{R}} }
\newcommand{\intRR}{\int\limits_{\mathbb{R}^2} }
\newcommand{\RR}{\mathbb{R}}
\newcommand{\xx}{{\mathbf{x}}}
\newcommand{\kk}{{\mathbf{k}}}
\newcommand{\XX}{\mathbf{X}}
\newcommand{\yy}{\mathbf{y}}
\newcommand{\YY}{\mathbf{Y}}
\newcommand{\zz}{\mathbf{z}}
\newcommand{\uu}{\mathbf{u}}
\newcommand{\xb}{\textbf{\textit{x}}} 
\newcommand{\yb}{\textbf{\textit{y}}} 
\newcommand{\zb}{\textbf{\textit{z}}}
\newcommand{\ub}{\textbf{\textit{u}}}
\newcommand{\ttheta}{\boldsymbol{\theta}}
\newcommand{\pphi}{\boldsymbol{\phi}}
\newcommand{\oomega}{\boldsymbol{\omega}}
\newcommand{\eeta}{\boldsymbol{\eta}}
\newcommand{\zzeta}{\boldsymbol{\zeta}}
\newcommand{\cchi}{\boldsymbol{\chi}}
\newcommand{\xxi}{{\boldsymbol{\xi}}}
\newcommand{\rmd}[1]{\mathrm d#1}
\newtheorem{thm}{Theorem}
\newtheorem{defi}[thm]{Definition}
\newtheorem{rmk}[thm]{Remark}
\newtheorem{cor}[thm]{Corollary}
\newtheorem{lem}[thm]{Lemma}
\newtheorem{prop}[thm]{Proposition}
\newtheorem{prob}{Problem}

\newcommand{\sr}[1]{\sout{\color{red}{#1}}}
\newcommand{\br}[1]{\textbf{\color{red} #1}}
\newcommand{\bb}[1]{\textbf{\color{blue} #1}}

\newtheorem{myproblem}{Problem}
\newtheorem*{mysubproblem}{Subproblem}
\newtheorem*{myhypothesis}{Hypothesis}

\newtheorem{assumption}{Theorem}

\newtheorem{assumptionalt}{Theorem}[assumption]

\newenvironment{assumptionp}[1]{
	\renewcommand\theassumptionalt{#1}
	\assumptionalt
}{\endassumptionalt}

\title{Implicit learning to determine variable sound speed and the reconstruction operator in photoacoustic tomography}
\author[1]{Gyeongha Hwang}
\author[2**]{Gihyeon  Jeon}
\author[3]{Sunghwan Moon}
\author[2*]{Dabin Park}
\affil[1]{Department  of  Mathematics,  Yeungnam  University,  Gyeongsan  38541,  Republic  of  Korea}
\affil[2]{School  of  Mathematics,  Kyungpook  National  University,  Daegu  41566,  Republic  of  Korea}
\affil[3]{Department  of  Mathematics,  Kyungpook  National  University,  Daegu  41566,  Republic  of  Korea}
\affil[*]{Corresponding  author: ekqls5264@knu.ac.kr}
\affil[**]{Corresponding  author: rydbr6709@knu.ac.kr}

\maketitle

\begin{abstract}
Photoacoustic tomography (PAT) is a hybrid medical imaging technique that offer high contrast and a high spatial resolution. One challenging mathematical problem associated with PAT is reconstructing the initial pressure of the wave equation from data collected at the specific surface where the detectors are positioned. The study addresses this problem when PAT is modeled by a wave equation with unknown sound speed $c$, which is a function of spatial variables, and under the assumption that both the Dirichlet and Neumann boundary values on the detector surface are measured. In practical, we introduce a novel implicit learning framework to simultaneously estimate the unknown $c$ and the reconstruction operator using only Dirichlet and Neumann boundary measurement data. The experimental results confirm the success of our proposed framework, demonstrating its ability to accurately estimate variable sound speed and the reconstruction operator in PAT.
\end{abstract}

Keywords: \texttt{photoacoustic tomography, unsupervised learning, inverse problem, wave equation}

\section{Introduction}

Photoacoustic tomography (PAT) is an imaging method that uses non-ionized laser pulses and ultrasound to produce detailed images of the internal structure of biological tissue.
The method is non-destructive, economical, and less harmful than other imaging options because it uses non-ionizing radiation \cite{Steinberg19}. For these reasons, it is used in a variety of biomedical applications, including skin melanoma detection, breast cancer detection,  blood oxygenation mapping, tumor angiogenesis monitoring, functional brain imaging, and methemoglobin measurement \cite{Xu06}.

In PAT, when the target object is irradiated with a non-ionizing laser pulse, the absorbed pulse generates a photoacoustic effect in which rapid thermal expansion leads to the production of acoustic waves, a phenomenon that was first discovered by Bell \cite{Bell80}.
The resulting acoustic waves are measured using ultrasound detectors as the data and used to reconstruct images of the target object.
This method is particularly useful for visualizing structures in optically opaque materials such as biological tissue, because it combines the high contrast of optical imaging with the high spatial resolution of ultrasound imaging (for more details, see \cite{Jiang18, Kuchment13, Xia14}).

One of the goals of PAT is to reconstruct initial pressure $f$, which potentially contains important biological information such as the presence and location of cancer cells, from the acoustic waves measured by the detector.
These acoustic waves, denoted by $p$, follow the wave equation:
\begin{equation}\label{eq:Wave}
	\left\{\begin{array}{ll}
		\partial_t^2 p(\xx,t) =  c(\xx) \triangle_\xx p(\xx,t) & (\xx,t) \in \RR^2 \times [0,\infty),\\
		p(\xx,t)|_{t = 0} = f(\xx) \text{ and } \partial_t p(\xx,t)|_{t=0} = 0 & \xx \in \RR^2
	\end{array}\right.
\end{equation}
Here, $c(\mathbf{x})$ represents the sound speed at location $\xx$, which is assumed to be continuous and bounded by two constants $c_{M}$ and $c_{m}$ such that $c_{M} \ge c(\mathbf{x}) \ge c_{m} \ge 0$. The measurements are made on the boundary of a region of interest $\Omega \subset \mathbb{R}^2$ under the reasonable assumption that $f$ has compact support within the bounded domain $\Omega$. We can then define wave forward operator $\mathcal{W}_c$, which maps initial pressure $f$ to solution $p$, i.e., $\mathcal{W}_c f = p$. The challenge then focuses on precisely reconstructing $f$ from the boundary measurements to derive an accurate internal image of the target object.

Most studies on the reconstruction of initial pressure $f$ from $\mathcal{W}_{c}{f}|_{\partial\Omega\times[0,\infty]}$ assume a known sound speed $c$.
The use of a constant $c$ has been studied analytically in \cite{Finch07, Xu05, Zangerl19} (for further mathematical details, see \cite{Ammari10,Kuchment08, Kuchment13} and references therein), while other analytical studies have employed a variable $c$.
Agranovsky and Kuchment study reconstruction of the initial pressure from Dirichlet data with a known variable sound speed in three-dimensional space \cite{Agranovsky07}.
Moon et al. also study the singular value decomposition of the wave forward operator with a known radial sound speed in $n$-dimensional space \cite{Moon23}.
Stefanov and Uhlmann address a more general problem with a Riemann metric instead of the sound speed \cite{Stefanov09}.
Other studies detailing numerical approaches include \cite{Belhachmi16, Hristova08, Qian11}.
Research has also been conducted on recovering $c$ when initial pressure $f$ is known \cite{Stefanov13}, while the sufficient conditions for simultaneously reconstructing $f$ and $c$ are discussed in \cite{Liu15}.

Recently, deep learning methods have been applied to PAT \cite{Suganyadevi22, Hauptmann20} image reconstruction \cite{Antholzer18}, handling limited data setups \cite{Antholzer19, Awasthi20, Gutta17, Jeon21, Shahid21}, and achieving a super-resolution \cite{Awasthi20}.
Many of these studies are based on supervised learning, for which initial pressure $f$ is the target data.
However, in practice, it cannot necessarily be assumed that the target data is known in PAT because the initial pressure represents the internal structure of the object.
Therefore, methods for training the reconstruction operator without using the target data need to be considered.
In line with this, in this paper, we propose an unsupervised learning method to implicitly estimate the variable speed of sound $c$ and the reconstruction operator.
The proposed method only utilizes paired Dirichlet and Neumann data boundary values and contributes to PAT research by demonstrating the use of implicit learning techniques.

The rest of this paper is structured as follows. The next section presents the formulation of our problem.
In Section \ref{sec:proposed method}, we describe our proposed framework and its loss function.
The numerical results are presented in Section \ref{sec:numerical simulations}.
Finally, Section \ref{section:conclusion} summarizes our research and its contributions to PAT, focusing on the variable speed of sound and initial pressure.

\section{Problem Formulation}
\label{subsection:problem setup}
Building on \eqref{eq:Wave}, we consider the problem of estimating variable sound speed $c$ and the reconstruction operator that recovers initial pressure $f$ from Dirichlet and Neumann boundary data.
We assume that $c \in C^{\infty}(\mathbb{R}^{2})$ satisfies $c_{M} \ge c(\mathbf{x}) \ge c_{m} \ge 0$.
Theorem 3.3 in \cite{Stefanov13} is a uniqueness theorem associated with this problem.
According to this theorem, for a given paired dataset of initial values and Dirichlet data, sound speed $c$ and the reconstruction operator can be uniquely determined (see the Appendix for further details of the problem formulation, methods, and experimental results for this type of paired dataset is given).
However, because it is an unreasonable assumption that the ground truth for the initial data can be used in PAT, obtaining paired data is impractical.
Therefore, in this paper, we address this problem without using initial data.

For convenience, Dirichlet data $\mathcal W_{c}f|_{\partial B \times [0,\infty)}$ is denoted as $\mathcal{D}_{c}f$ and Neumann data $\partial_\nu \mathcal W_{c}f|_{\partial B \times [0,\infty)}$ is denoted as $\mathcal{N}_{c}f$, 
where $\nu$ is the unit outward normal vector at the detector surface $\partial B$ and $B$ is the unit ball in $\mathbb{R}^{2}$.

\begin{myproblem}
	When a collection of Dirichlet and Neumann data pairs
	\begin{equation}\label{eq:data}
		\left\{ (\mathcal{D}_{c}f, \mathcal{N}_{c}f): f \in L^{2}(\RR^2) \text{ with } \operatorname{supp} f \subset B \right\} 
	\end{equation}
	is given, estimate the sound speed $c$ and the reconstruction operator $\mathcal{D}_c^{-1}$. 
\end{myproblem}

When $c$ is known, {\bf Problem 1} aligns with the uniqueness theorem ({\bf Theorem A}):
\begin{assumptionp}{A}(\cite[Theorem 8]{Agranovsky07})
	For a known sound speed $c$, 
	the initial pressure $f$ is uniquely determined by the Dirichlet boundary value $\mathcal D_c f$.
\end{assumptionp}

We can now discuss the uniqueness of $c$ in \textbf{Problem 1}. For this, we recall the Calder\'{o}n problem on the conductivity equation \cite{Feldman19, Sylvester90}: For the equation  
\begin{equation}\label{eq:conductivity}
		\nabla \cdot (\gamma(\xx)\nabla u(\xx)) = 0 \qquad \xx \in B,
\end{equation}
with Dirichlet condition $u|_{\partial B} = g$, the problem is 
uniquely determining conductivity function $\gamma$ from the knowledge of the (bounded linear) Dirichlet to Neumann map 
$\Lambda_{\gamma} : H^{1/2}(\partial B) \to H^{-1/2}(\partial B)$ defined by
\begin{equation*}
	\Lambda_{\gamma}(g) = \partial_\nu u|_{\partial B},
\end{equation*}
where $H^\alpha(B), \alpha \in \mathbb{R}$ is the Sobolev space.

Note that conductivity equation \eqref{eq:conductivity} is equivalent to 
\begin{equation}\label{eq:Schrodinger}
		(-\Delta + q(\xx))w(\xx) = 0 \qquad\xx  \in B,
\end{equation}
where $q(\xx) = \dfrac{\nabla \gamma^{1/2}(\xx)}{\gamma^{1/2}(\xx)}$ and $w(\xx) = \gamma^{1/2}(\xx)u(\xx)$.

Based on wave equation \eqref{eq:Wave}, we propose the following conjecture.

\textbf{Conjecture:}
Assume that $c_{1}$ and $c_{2}$ are the known constant $c_0$ on $B^c$. 
If, for any $f\in L^2(B)$, 
\begin{equation*}
	\mathcal{D}_{c_1} f = \mathcal{D}_{c_2} f \text{ and } \mathcal{N}_{c_1} f = \mathcal{N}_{c_2} f,
\end{equation*}
then $c_1 = c_2$ on $\RR^2$.

Because $c_{1}$ and $c_{2}$ are known constant $c_0$ on $B^c$, it suffices to show that $c_1 = c_2$ on $B$.
Let us define 
\begin{equation*}
\overline{p}(\xx,t) := 
\left\{ \begin{array}{ll} p(\xx,t)|_{\overline{B} \times [0,\infty)} & t \geq 0,\\
p(\xx,-t)|_{\overline{B} \times [0,\infty)} & t < 0,
\end{array}\right.
\end{equation*}
where $p$ is the solution for \eqref{eq:Wave}.
Then, $\overline{p}$ satisfies 
\begin{equation}\label{eq:modwave}
	\left\{\begin{array}{ll}
		\partial_t^2 \overline{p}(\xx,t) =  c(\xx) \Delta_\xx \overline{p}(\xx,t) & (\xx,t) \in B \times \RR,\\
		\overline{p}(\xx,0) = f(\xx) \text{ and } \partial_t \overline{p}(\xx,t)|_{t=0} = 0 & \xx \in B.
	\end{array}\right.
\end{equation}
Taking the Fourier transform of \eqref{eq:modwave} with respect to the time variable $t$, we have
\begin{equation*}
	\left(-\Delta_\xx + \frac{(2\pi \text{i} \omega)^2}{c(\xx)}\right)\mathcal{F}_t(\overline{p}(\xx,\cdot))(\omega) = 0. \nonumber
\end{equation*}
The Calderón problem appears to give us the uniqueness of $c$ in \eqref{eq:modwave} given \eqref{eq:data}. 
However, there are gaps in this application. For $f\in  L^2(B)$ we cannot be sure that $\mathcal{D}_c f \in H^{1/2} (\partial B)$.
Moreover, it is not easy for set $\{ \mathcal{D}_c f : f \in L^2(B) \}$ to be equal to $H^{1/2} (\partial B)$.
Despite these gaps, our conjecture is validated  experimentally in Section \ref{sec:numerical simulations}.

In practical terms, we focus on the following problem, emphasizing that only a finite collection of data can be used in real-world applications.

\begin{myproblem}
	When a sufficiently large finite collection of Dirichlet and Neumann data pairs
	\begin{equation}\label{eq:finitedata}
		\left\{ (\mathcal{D}_{c}f_{i}, \mathcal{N}_{c}f_{i})_{i = 1, \dots, N}: f_{i} \in L^{2}(\RR^2) \text{ with } \operatorname{supp} f_{i} \subset B \right\} 
	\end{equation}
	is given, estimate the sound speed $c$ and the reconstruction operator $\mathcal{D}_c^{-1}$. 
\end{myproblem}

Our problem has the following difficulties:
\begin{enumerate}
	\item Target data $f$ in PAT is unavailable, and
	\item The explicit formula of the wave forward operator is  unknown.
\end{enumerate}

In this paper, we propose an implicit learning method for estimating $c$ and $\mathcal{D}_c^{-1}$.
The proposed method uses a paired dataset of Dirichlet and Neumann data and an iterative method for the wave forward operator.

\section{Proposed method}
\label{sec:proposed method}

\begin{figure}[!h]
	\centering
    \mbox{\includegraphics[width=1\textwidth]{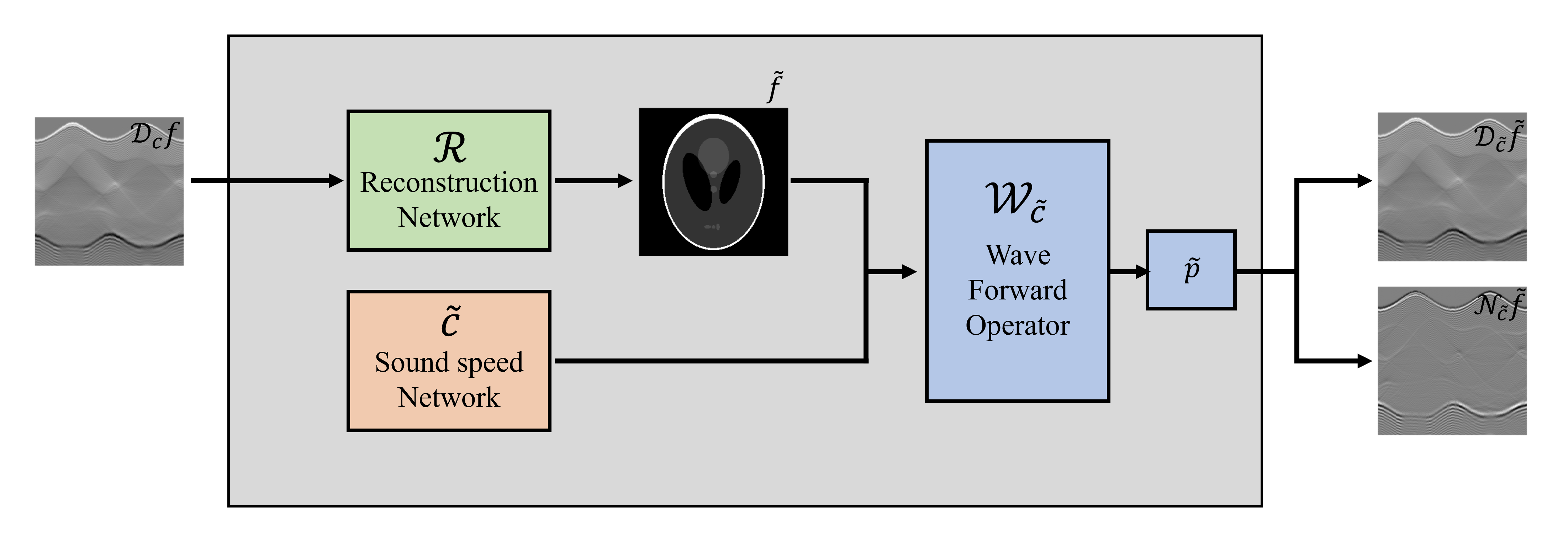}}
	\caption{Proposed framework}
	\label{fig:framework}
\end{figure}

Our goal is to estimate $c$ and $\mathcal{D}_c^{-1}$ simultaneously from given data set $\mathscr{T} = \{(\mathcal{D}_{c}f_{i}, \mathcal{N}_{c}f_{i})\}_{i = 1}^{N}$.
The proposed framework is shown in Figure \ref{fig:framework} and consists of three components:

\begin{enumerate}
	\item[1.] Sound speed network $\tilde{c}$
	\item[2.] Reconstruction network $\mathcal R$
	\item[3.] Wave forward operator $\mathcal W_{\tilde{c}}$
\end{enumerate}

Sound speed network $\tilde{c}$  is designed to estimate the unknown $c$ in \eqref{eq:Wave}. Simultaneously, reconstruction network $\mathcal{R}$ approximates ${\mathcal{D}_{c}^{-1}}$. This network takes Dirichlet data $\mathcal{D}_c f$ and outputs a reconstruction of initial data $\tilde{f}$. Wave forward operator $\mathcal{W}_{\tilde{c}}$, using an iterative scheme, solves wave equation \eqref{eq:Wave} with the estimated $\tilde{c}$ and the reconstructed $\tilde{f}$. From solution $\tilde p$, Dirichlet $\mathcal D_{\tilde c}\tilde f$ and Neumann $\mathcal N_{\tilde c}\tilde f$ data are obtained.

If the proposed conjecture is valid, then 
\begin{equation}
	||\mathcal{D}_{\tilde{c}} \tilde{f} - \mathcal{D}_{c} f||_2 = 0 \quad \text{and} \quad ||\mathcal{N}_{\tilde{c}} \tilde{f} - \mathcal{N}_{c} f||_2 = 0 \text{ for all } f \in L^2(B) \nonumber
\end{equation} 
implying that $\tilde{c}$ and $c$ are the same.
{\bf Theorem A} guarantees that for $\tilde{c} = c$, if $||\mathcal{D}_{\tilde{c}} \tilde{f} - \mathcal{D}_{c} f||_2 = 0$ holds, then $\tilde{f} = f$.  
Thus, for given data set $\mathscr{T} = \{(\mathcal{D}_{c}f_{i}, \mathcal{N}_{c}f_{i})\}_{i = 1}^{N}$, we define the following loss function:
\begin{equation}\label{eq: loss}
	\mathcal{L}_{A}(\mathscr{T}) = \sum \Big[ \lambda_{\mathcal{D}}|| \mathcal{D}_{\tilde{c}}\tilde{f}_{i} -\mathcal{D}_{c}f_{i}||_{2}^{2} +  \lambda_{\mathcal{N}} || \mathcal{N}_{\tilde{c}}\tilde{f}_{i}-\mathcal{N}_{c}f_{i}||_2^2 \Big].
\end{equation}
To effectively handle noisy data, we also incorporate a Total Variation (TV) regularization term for the reconstructed image \cite{Rudin92}:
$$
\mathcal{L}_{B}(\mathscr{T}) = \lambda_{\text{TV}} \sum ||\nabla \tilde{f}_{i}||_{1} \qquad \text{where } \tilde{f}_{i} = \mathcal{R}(\mathcal{D}_{c}f_{i}).
$$
The overall loss function is then given by:
$$
\mathcal{L}(\mathscr{T}) = \mathcal{L}_{A}(\mathscr{T}) + \mathcal{L}_{B}(\mathscr{T}).
$$
In the following subsection, we describe each component of the architecture in detail.

\subsection{Sound speed network} 
\label{subsection:sound speed network}
We propose a neural network to estimate the unknown sound speed $c$ in \eqref{eq:Wave}.
This strategy is based on the universal approximation theorem \cite{Cybenko89}, which asserts that neural networks can approximate any continuous function.
For this task, we employ a Multilayer Perceptron (MLP) architecture.
Our MLP is structured to take two-dimensional spatial input $\mathbf{x}$ and output the estimated sound speed $\tilde{c}(\mathbf{x})$.
The MLP consists of three hidden layers, each containing 50 hidden units.
We choose the sine function as the activation function for the hidden layers due to its periodic nature, which enhances the model's ability to accurately and efficiently represent complex natural signals and their derivatives \cite{Sitzmann20}.
Considering that the sound speed has known upper and lower bounds, we use a hyperbolic tangent function to constrain the output, enabling it to fit within the predetermined speed range.
Additionally, we set the sound speed values outside of the detector's range to a known constant.

\subsection{Reconstruction network}
The reconstruction network is designed to approximates the reconstruction operator that reconstructs initial data $f$ using Dirichlet data $\mathcal{D}_c f$.
Because $\mathcal{D}_c$ is linear, its inverse is also linear, thus the reconstruction operator can be approximated by a neural network with a single linear layer.
However, when the input has a high-resolution, this structure leads to a sharp increase in the number of neural network parameters, which grows quadratically with the resolution of input image.
This results in substantial memory demands and longer computation times.
For example, a low-resolution image of size $64\times64$ requires 16,777,216 parameters; conversely, a high-resolution image of size $256\times256$ requires 4,294,967,296 parameters, a 256-fold increase.
An excessive number of parameters can also increase the risk of overfitting, adversely affecting the generalization capability of the model.

To mitigate these challenges, our proposed network architecture incorporates techniques such as downsampling and upsampling to significantly reduce the number of parameters while maintaining the network's efficacy.
This approach optimizes both the computational efficiency and memory usage, allowing it to effectively handle higher-resolution data.
The architecture of the reconstruction network is detailed in following Figure \ref{fig:R_network}.

\begin{figure}[h]
	\centering
    \mbox{\includegraphics[width=1\textwidth] {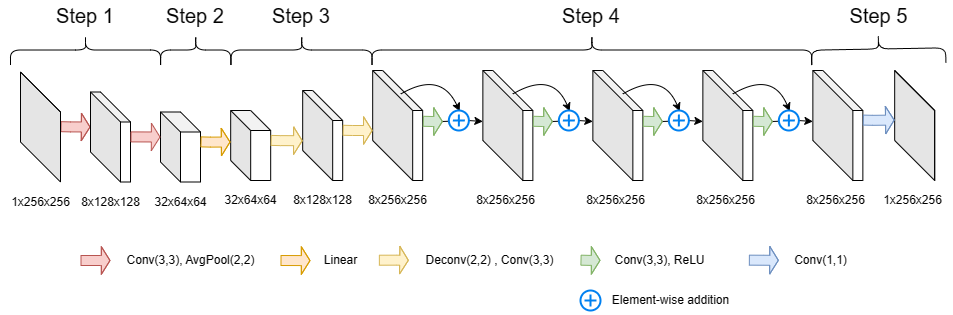}}
	\caption{\textbf{Reconstruction network $\mathcal R$.}
	The Conv(a, b) layer is a 2D convolution layer with a kernel size of (a, b), and a stride of (1, 1) while padding is employed to ensure that the output is the same size as the input.
	The AvgPool(a, b) layer is a 2D average pooling layer with a kernel size of (a, b).
	The linear layer is a fully connected layer that transforms an input with dimension of $64\times64$ to an output with the same dimensions.
	The Deconv(a, b) layer is a 2D transposed convolution layer with a kernel size of (a, b) and a stride of (2, 2).
	Each layer contains no bias.
	}
	\label{fig:R_network}
\end{figure}

The core process for the reconstruction network is as follows:

\begin{enumerate}
\item[\textbf{Step 1.}] 
						The input data (sinogram) undergoes a reduction in resolution via the convolution and average pooling layers, converting it into low-resolution, multi-channel data.
\item[\textbf{Step 2.}] 
						The low-resolution data from Step 1 is fed through a fully connected layer, which produces further low-resolution, multi-channel image domain data.
						This step relies on the linearity of the reconstruction operator.
\item[\textbf{Step 3.}] 
						Deconvolution and convolution layers are employed to upscale the data from Step 2, resulting in high-resolution, multi-channel image domain data.
						This step primarily focuses on effectively enhancing the resolution of the final image output.
\item[\textbf{Step 4.}] A residual network (Res-Net) is used to further refine the reconstruction quality.
						Res-Net can significantly enhance the fidelity of reconstructed images.
\item[\textbf{Step 5.}] 
						To ensure that final image $f$ has values within the physiological range of $[0, 1]$, hyperbolic tangent transformation is employed.
						This function $x \mapsto 0.5 \tanh(x - 0.5) + 0.5$ helps to calibrate the output to the desired range via normalization.
\end{enumerate}

\subsection{Wave forward operator}

Because the explicit form of the wave forward operator is unknown, we use an iterative method.
In particular, to calculate the wave propagation, we use the well-known $k$-space method \cite{Cox05, Cox07}.
The approximation of wave propagation for the subsequent time step is obtained from equation \cite{Hwang23}:
\begin{equation}\label{eq: wave propagation}
p(\xx,t+\triangle t) = 2 p(\xx,t)  -  p(\xx,t-\triangle t) - \tilde{c}(\xx) \mathcal F_\kk^{-1}\left[4\sin^2\left(\dfrac{(\triangle t)|\cdot|}2\right)  \mathcal F_\xx[p](\cdot,t)\right](\xx),
\end{equation}
where $\mathcal F_\xx$ is the Fourier transform, and $\mathcal F_\kk^{-1}$ is the inverse Fourier transform.
In this formulation, $p(\xx,t+\triangle t)$ denotes the solution for the subsequent time step, illustrating how the system evolves over time, particularly under the influence of the spatially dependent term $\tilde{c}(\xx)$.

\section{Numerical Results}
\label{sec:numerical simulations}

We use Shepp-Logan phantoms for our simulation data. 
A Shepp-Logan phantom, introduced by Shepp and Logan in 1974 \cite{Shepp74}, is an artificial image representing a cross-section of the brain.
It consists of several ellipses each defined by six parameters: 
the center coordinates of the ellipse, the major axis, the minor axis, the rotation angle, and the intensity. 
We create a set of phantoms $\left\{ F_{i} \right\}_{i=1}^{4,096}$ by changing these parameters and obtain the Dirichlet data set 
\begin{equation}
	\left\{ D_{c}^{r}F_{i} : \text{Dirichlet data for initial $F_{i}$ and sound speed $c$ on a ball of radius $r$ }\right\}_{i=1}^{4,096} \nonumber
\end{equation}
for a given sound speed $c$ by applying the iterative wave propagation method \eqref{eq: wave propagation} to the phantoms. 
Noisy data are generated by adding Gaussian noise to 1\% of the maximum value of the original data. 
We use 2,048 phantoms as a training set, 1,024 as a validation set, and the remaining 1,024 as a test set.

We set  $c(x, y; h) = \exp(-h(x^2 + y^2))$. Simulations are conducted using two different sound speed profiles:
\begin{equation}\label{def: sound speed}
	\begin{aligned}	
		\text{Type 1: }& c_{1}(x, y) = \frac{3}{4}\Big( c(x - 0.5, y - 0.6; 2) + c(x - 0.5, y + 0.5; 2) + c(x + 0.6, y - 0.45; 2) \\
		& \qquad\qquad\quad + c(x + 0.55, y + 0.6; 2) \Big), \\
		\text{Type 2: }& c_{2}(x, y) = \frac{23}{35}\Big( c(x , y - 0.675; 2.5) + c(x - 0.725, y - 0.175; 2.5) + c(x - 0.4, y + 0.625; 2.5) \\ 
		& \qquad\qquad\quad\quad + c(x + 0.625, y - 0.225; 2.5) + c(x + 0.45, y + 0.575; 2.5)\Big),
	\end{aligned}
\end{equation}
which are characterized by a constant value outside of $B$.

In practice, for small a $h > 0$, the Neumann data at the boundary of $B$ can be approximated by
\begin{equation}\label{eq: neumann_approx}
	\partial_{\boldsymbol{\nu}_{\xx}} \mathcal{W}_{c}f (\xx, t) \approx \frac{\mathcal{W}_{c}f(\xx, t) - \mathcal{W}_{c}f(\xx - h \boldsymbol{\nu}_{\xx}, t)}{h},\quad \xx \in \partial B
\end{equation}
where $\boldsymbol{\nu}_{\xx}$ is the unit outward normal vector at $\xx \in \partial B$. Having Dirichlet and Neumann data pairs is equivalent to having Dirichlet data pairs ${(D_{c}^1 F_{i}, D_{c}^{1-h} F_{i})}$. Therefore, we modify the loss function described in \eqref{eq: loss} as follows:
\begin{equation}\label{eq: actual loss}
	\mathcal{L} = \lambda_{D}\sum\limits_{i} \left\| D_{c}^{1}F_{i}- D_{\tilde{c}}^{1}\tilde{F}_{i} \right\|_{2}^{2} + \lambda_{N} \sum\limits_{i}\left\| D_{c}^{1-h}F_{i}- D_{\tilde{c}}^{1-h}\tilde{F}_{i} \right\|_{2}^{2} + \lambda_{\text{TV}} \sum\limits_{i}||\tilde{F}_{i}||_{\text{TV}},
\end{equation}
where $\tilde{F} = \mathcal{R}(D_{c}^{1}F)$.

We choose $h = 0.05$, $\lambda_{D} = 1$, $\lambda_{N} = 1$, and $\lambda_{\text{TV}} = \sigma/5$, where $\sigma$ is the noise level.
To minimize the loss function \eqref{eq: actual loss}, we use the Adam optimizer \cite{Kingma14} with a learning rate of $0.001$ and a batch size of $2$.
The training takes $10^{5}$ iterations. 

Table \ref{table: error} shows the relative errors of the reconstructed image and the approximated sound speed. 
The results show that the proposed framework effectively reconstructs both the initial pressure and sound speed with high accuracy, maintaining its performance even with noisy data.

\begin{table}[H]
	\caption{Relative error of the reconstructed image and sound speed}
	\label{table: error}
	\setlength{\tabcolsep}{10pt}
	\renewcommand{\arraystretch}{1.5}
	\centering
	\begin{tabular}{c|c|c|c|c}
		\hline
		\multirow{2}{*}{Sound speed profile} & \multicolumn{2}{c}{Noise-free data} & \multicolumn{2}{c}{Noisy data} \\
		\cline{2-5}
		& $\textstyle \frac{1}{N}\sum\| f_{i} \|_{\text{rel}}$ & $\| c \|_{\text{rel}}$ & $\textstyle \frac{1}{N}\sum\| f_{i} \|_{\text{rel}}$ & $\| c \|_{\text{rel}}$ \\
		\hline
		\hline
		Type 1 & 0.0621 & 0.0014 & 0.0785 & 0.0014 \\
		Type 2 & 0.0613 & 0.0017 & 0.0759 & 0.0028 \\
		\hline
	\end{tabular}
\end{table}

\begin{figure}[H]
    \centering
	\includegraphics[width=0.7\textwidth ]{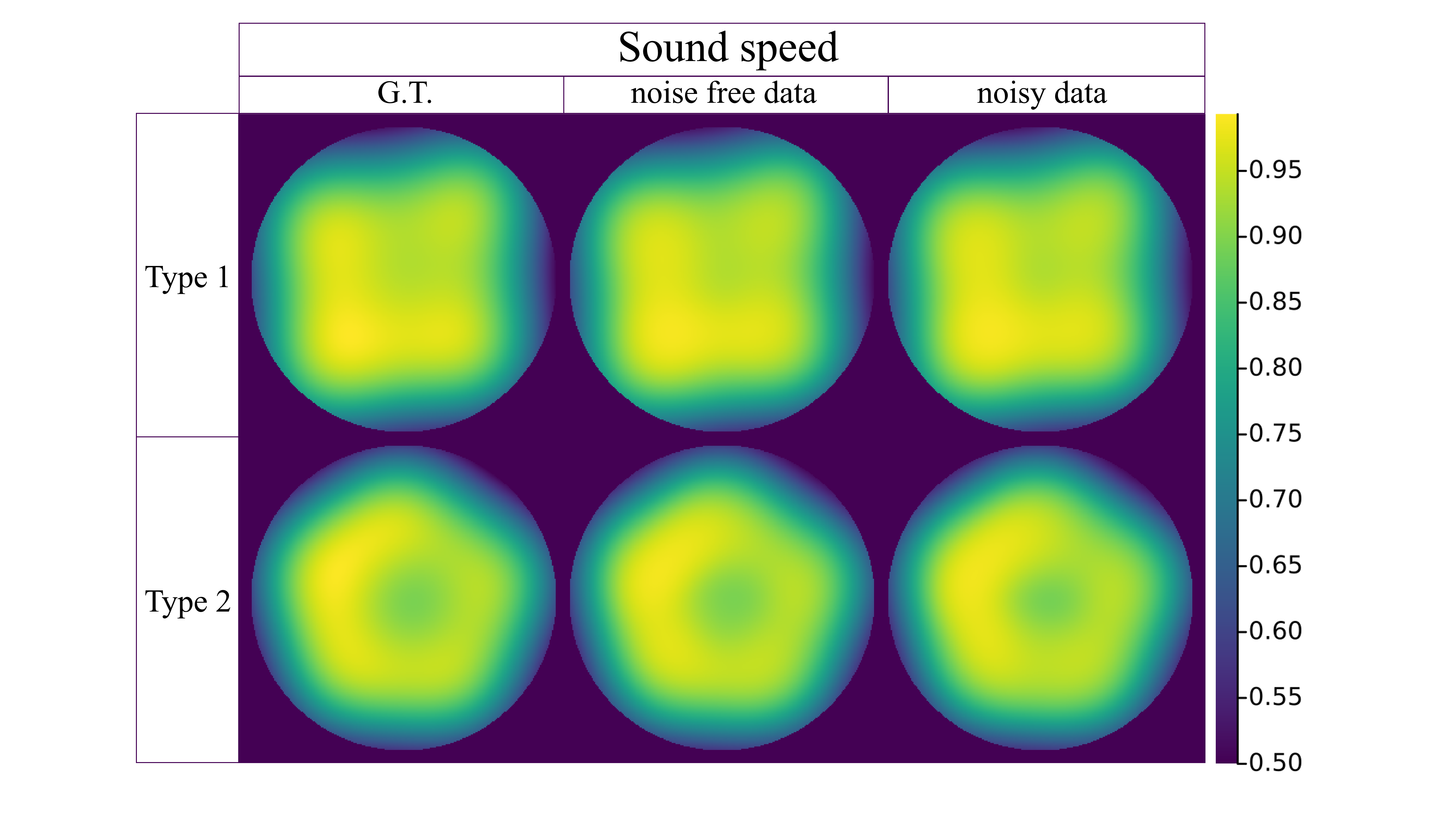}
    \caption{Estimated results for Type 1 and Type 2 sound speed profiles. The first, second, and third columns represent the ground truth, the estimated speed for noise-free data, and the estimated speed for noisy data, respectively.}
    \label{fig:estimated_c}
\end{figure}

\begin{figure}[H]
	\centering
	\subfloat[Results for Type 1 speed profile]{\includegraphics[width=0.7\textwidth ]{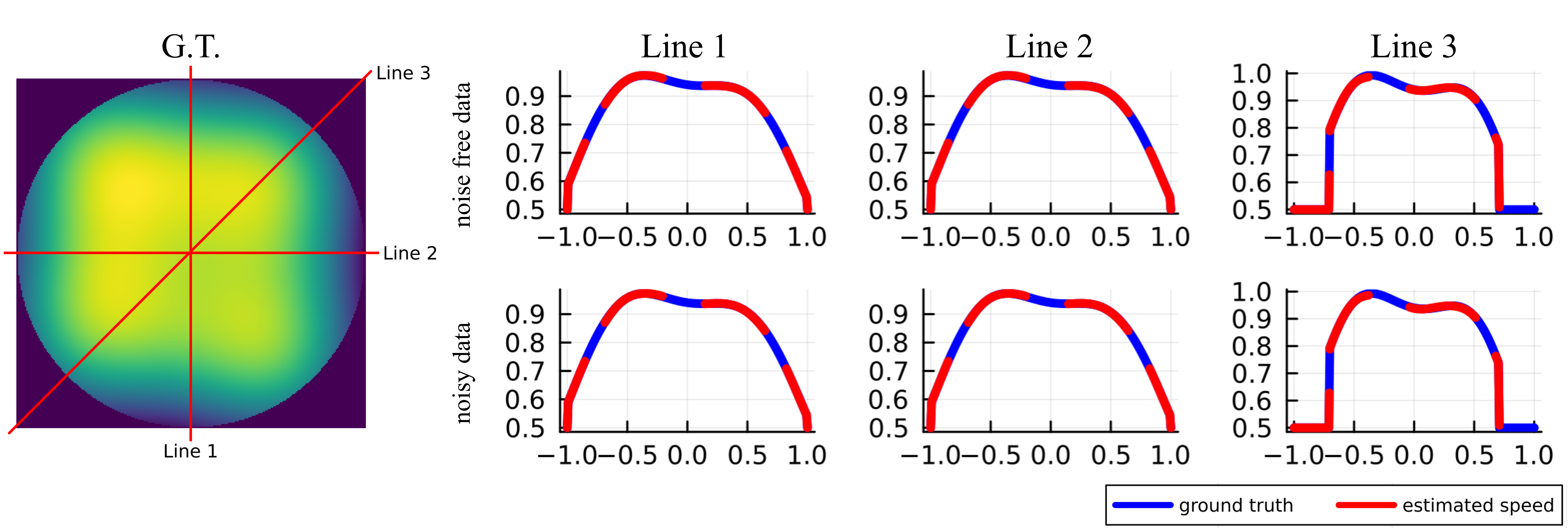}} \\
	\subfloat[Results for Type 2 speed profile]{\includegraphics[width=0.7\textwidth ]{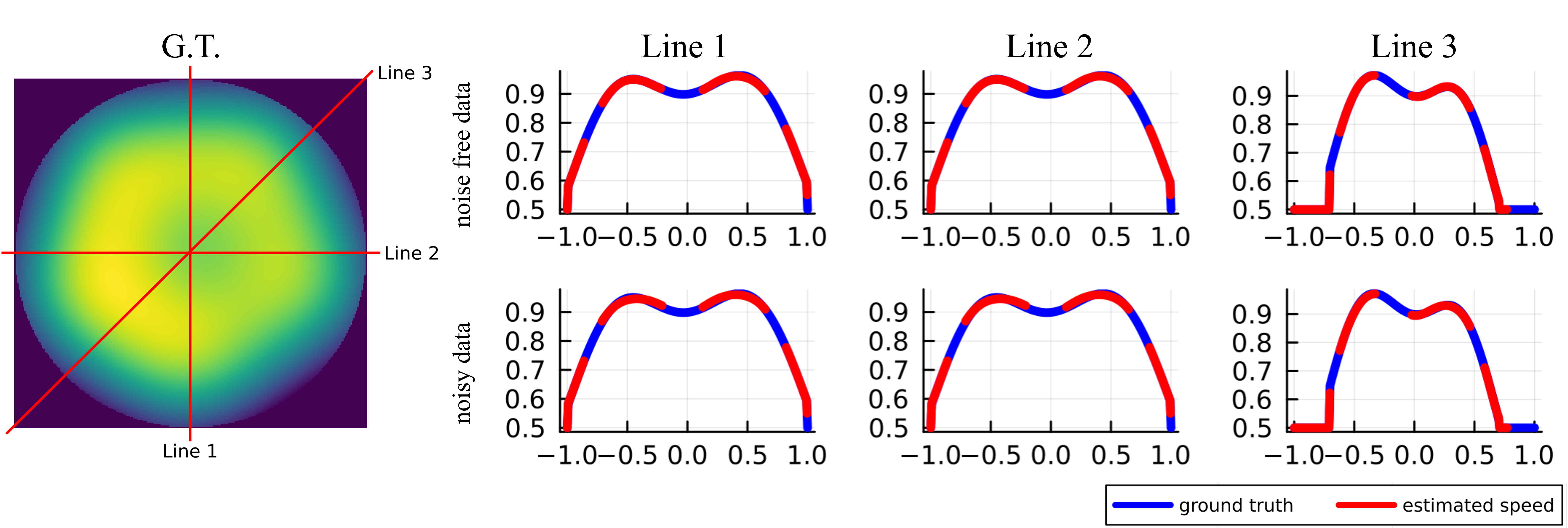}}
    \caption{In the left image, the ground truth is shown with three lines labeled Line 1, Line 2, and Line 3. The right images present cross-sectional views corresponding to these lines. The blue lines represent the ground truth, while the red dashed lines represent the estimated speed.}
    \label{fig:sliced_c}
\end{figure}

Figures \ref{fig:estimated_c} and \ref{fig:sliced_c} show the ground truth of the speed of sound along with their approximation output from the sound speed network. The relative error for $f$ with respect to ground truth $f_{\text{gt}}$ is defined by $\| f \|_{\text{rel}} = \dfrac{\| f_{\text{gt}} - f \|_{2}}{\| f_{\text{gt}} \|_{2}}$. The results demonstrate the ability of the network to accurately approximate the sound speed.
This successful approximation over different sound speed scenarios ($c_1, c_2$) illustrates the adaptability and robustness of the proposed network under varying conditions.

\begin{figure}[H]
	\centering
	\includegraphics[width=1\textwidth]{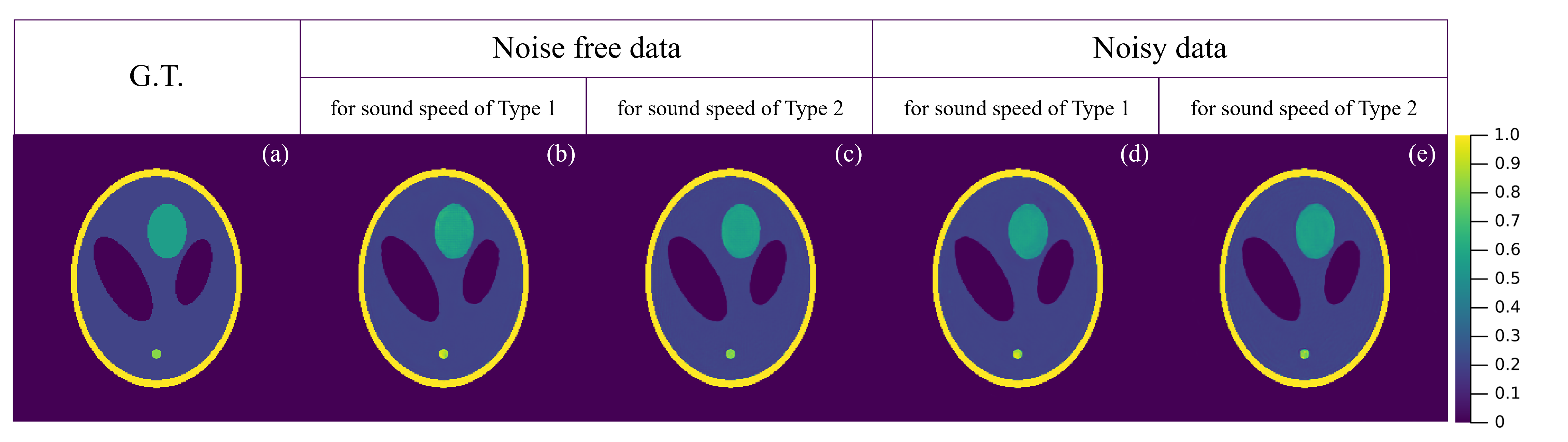}
	\caption{Reconstruction results: (a) ground truth (GT);
	(b) and (c) output of the reconstruction network for noise-free data for Type 1 and Type 2 sound speed profiles;
	(d) and (e) output of the reconstruction network for noisy data for Type 1 and Type 2 sound speed profiles.}
	\label{fig: reconstructed f}
\end{figure}

Figure \ref{fig: reconstructed f} illustrates the ground truth and reconstruction results for the reconstruction network, for both noise-free and noisy data sets. 

\section{Conclusion}
\label{section:conclusion}

In this study, we introduce a framework for addressing the problem in PAT associated with the simultaneous estimation of sound speed $c$ and reconstruction operator $\mathcal{D}_c^{-1}$. Our approach employs implicit learning to approximate both $c$ and $\mathcal{D}_c^{-1}$ using only Dirichlet and Neumann boundary data. This greatly reduces the dependence on large labeled datasets that is typical of medical imaging.

\section*{Appendix}
\label{sec:appendix}
In this Appendix, we address the problem of estimating the sound speed $c \in C^{\infty}(\mathbb{R}^{2})$ and the reconstruction operator $\mathcal{D}_{c}^{-1}$ when a paired dataset of initial pressure $f$ and Dirichlet data $\mathcal{D}_{c}f$ is given. With access to a sufficiently large dataset $\left\{ (f_{i}, \mathcal{D}_{c}f_{i}) \right\}_{i=1}^N$, we can estimate the reconstruction operator using supervised learning. This involves minimizing the following loss function:
$$
\frac{1}{N} \sum_{i} \| f_{i} - \mathcal{R}(\mathcal{D}_{c}f_{i}) \|_{2}^{2}.
$$
Here, $\mathcal{R} : \mathcal{D}_{c}f \mapsto f$ represents a neural network designed to approximate the reconstruction operator $\mathcal{D}_{c}^{-1}$. 

To estimate the sound speed $c(\cdot)$, we recall the following uniqueness theorem.
\begin{assumptionp}{B}(\cite[Theorem 3.3]{Stefanov13})
Let $c$ and $\tilde{c}$ be two smooth positive speeds equal to 1 outside $\Omega$. Let $\sum_{s}$, $s_{1} \le s \le s_{2}$ be a continuous family of compact oriented surfaces\footnote[1]{For a detailed definition of $\sum_{s}$, see \cite[pages 9--11]{Stefanov13}}. Let
$$
\mathcal{D}_{c} f = \mathcal{D}_{\tilde{c}}f \quad on \quad [0, T] \times \partial \Omega, \qquad with\ T > \max\limits_{s}\operatorname{dist} ( {\textstyle \sum_{s}} \cap \overline{\Omega}, \partial \Omega).
$$
Assume that for some compact $K \subset \overline{\Omega}$,
$$
\operatorname{supp} (\tilde{c} - c) \subset K, \qquad \Delta f(x) \ne 0 \quad for\ x \in K.
$$
Then $\tilde{c} = c$ in $\bigcup \sum_{s}$. If, in particular, $\bigcup \sum_{s}$ is dense in $K$, and $T > \operatorname{dist} (\Omega, \overline{\Omega})$, then $\tilde{c} = c$.
\end{assumptionp}
From \cite[Theorem 3.3]{Stefanov13}, we have:
\begin{equation*}
	\text{If } \mathcal{D}_{\tilde{c}}f = \mathcal{D}_{c}f \text{ for all }f \in L^2(B), \text{ then } \tilde{c} = c \text{ on } B.
\end{equation*} 
Based on this discussion, we define the loss function as follows:
\begin{equation}\label{eq:loss_appendix}
	\mathcal{L} = \lambda_{\mathcal{D}}||\mathcal{D}_{\tilde{c}}\tilde{f} - \mathcal{D}_{c}f||_2^2 + \lambda_{I}||\tilde{f} - f||_2^2 + \lambda_{\text{TV}}||\tilde{f}||_{\text{TV}},
\end{equation} 
where $\tilde{f} = \mathcal{R}(\mathcal{D}_{c}f)$. 

\begin{figure}[H]
    \centering
	\includegraphics[width=0.7\textwidth ]{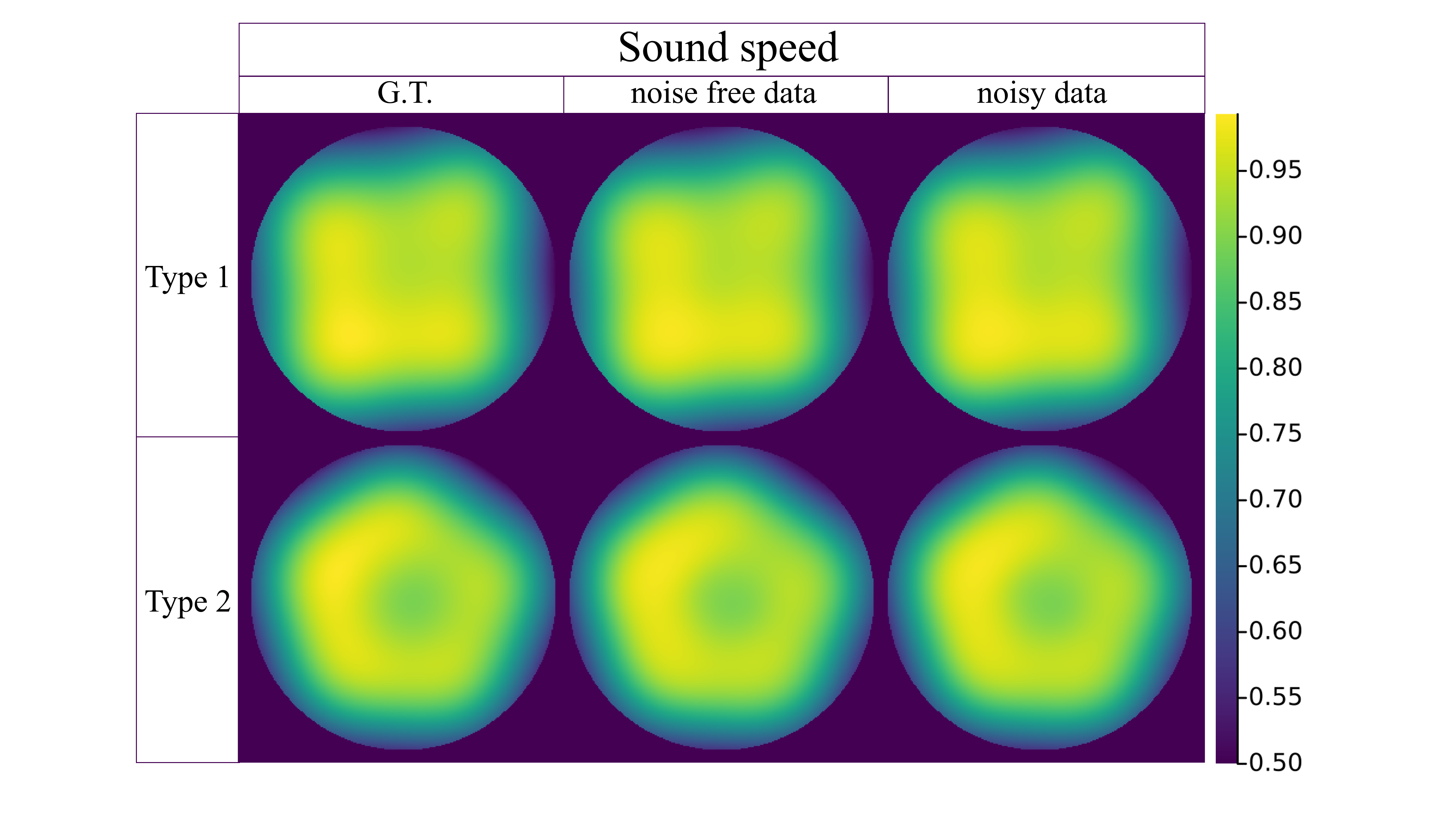}
    \caption{Estimated results for Type 1 and Type 2 sound speed profiles. The first, second, and third columns present the ground truth, the estimated speed for noise-free data, and the estimated speed for noisy data, respectively.}
    \label{fig:estimated_c_appendix}
\end{figure}

\begin{figure}[H]
	\centering
	\includegraphics[width=1\textwidth]{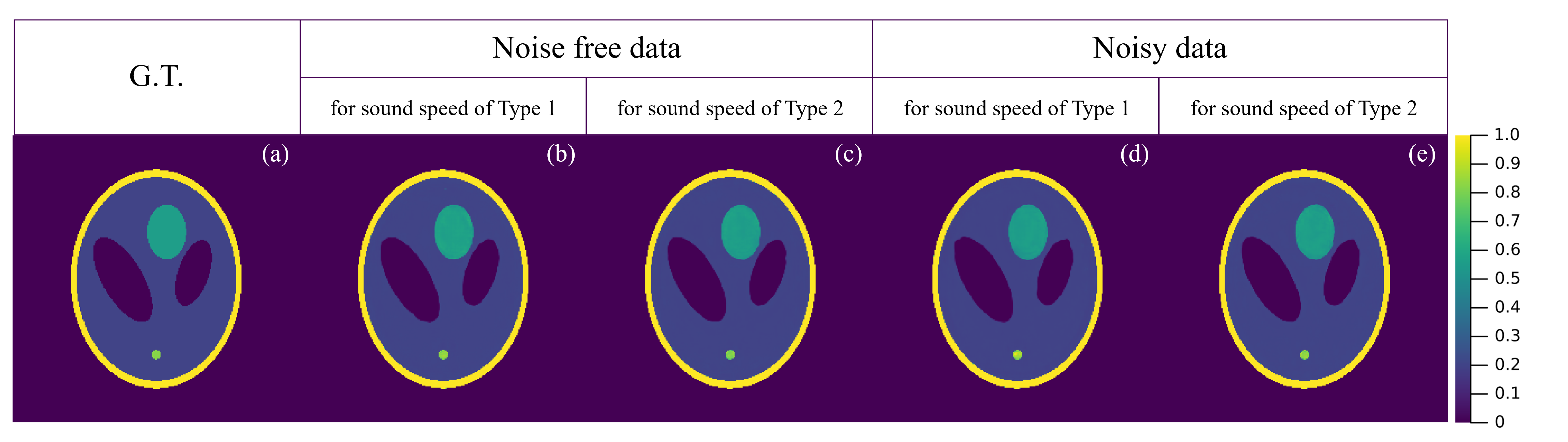}
	\caption{Reconstruction results: (a) ground truth (GT);
	(b) and (c) output of the reconstruction network for noise-free data for Type 1 and Type 2 sound speed profiles;
	(d) and (e) output of the reconstruction network for noisy data for Type 1 and Type 2 sound speed profiles.}
	\label{fig: reconstructed f_appendix}
\end{figure}

Experiments are conducted using the same dataset and framework as described in the main text. For detailed results, please refer to Table \ref{table: error_appendix}, Figure \ref{fig:estimated_c_appendix}, and Figure \ref{fig: reconstructed f_appendix}. The numerical results in the main text and Appendix show that the proposed implicit learning method, which does not use explicit target data, performs similarly to supervised learning. Additionally, the results in the Appendix are derived by simply modifying the loss function within the framework proposed in the main text, further confirming the applicability and robustness of our approach.

\begin{table}[H]
	\caption{Relative error of the reconstructed sound speed} 
	\label{table: error_appendix}
	\setlength{\tabcolsep}{10pt}
	\renewcommand{\arraystretch}{1.5}
	\centering
	\begin{tabular}{c|c|c|c|c}
		\hline
		\multirow{2}{*}{Sound speed profile} & \multicolumn{2}{c}{Noise-free data} & \multicolumn{2}{c}{Noisy data} \\
		\cline{2-5}
		& $\textstyle \frac{1}{N}\sum\| f_{i} \|_{\text{rel}}$ & $\| c \|_{\text{rel}}$ & $\textstyle \frac{1}{N}\sum\| f_{i} \|_{\text{rel}}$ & $\| c \|_{\text{rel}}$ \\
		\hline
		\hline
		Type 1 & 0.0605 & 0.0016 & 0.0766 & 0.0015 \\
		Type 2 & 0.0572 & 0.0020 & 0.0767 & 0.0019 \\
		\hline
	\end{tabular}
\end{table}

\begin{figure}[H]
	\centering
	\subfloat[Results for Type 1 speed]{\includegraphics[width=0.7\textwidth ]{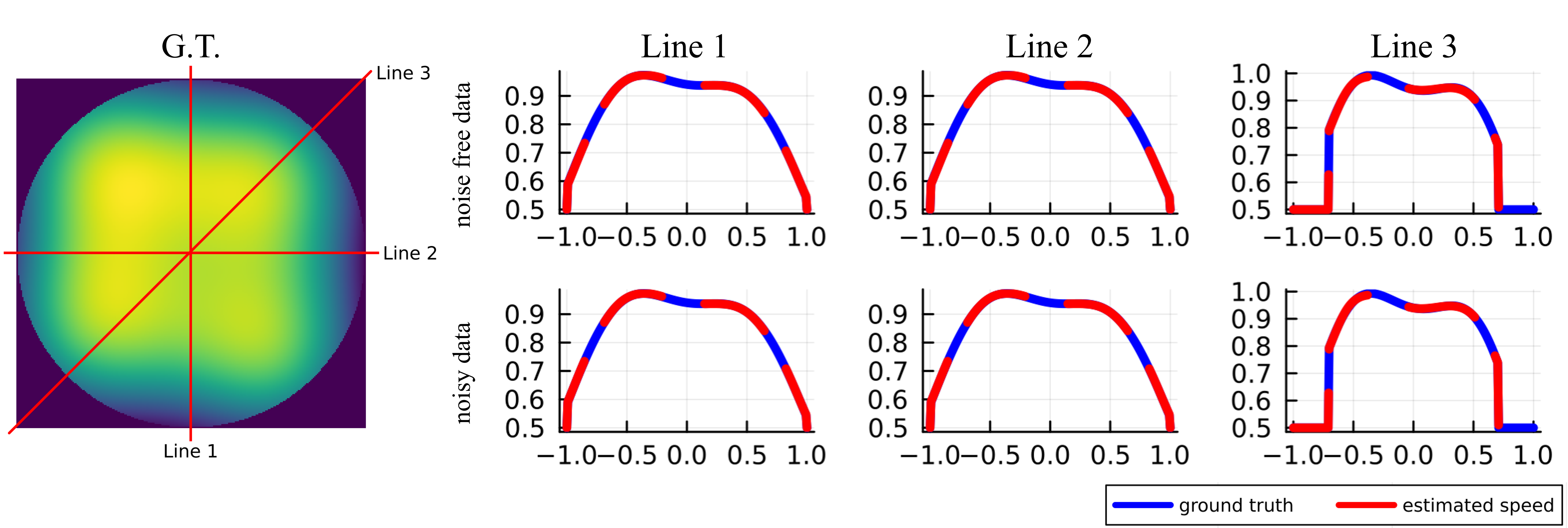}} \\
	\subfloat[Results for Type 2 speed]{\includegraphics[width=0.7\textwidth ]{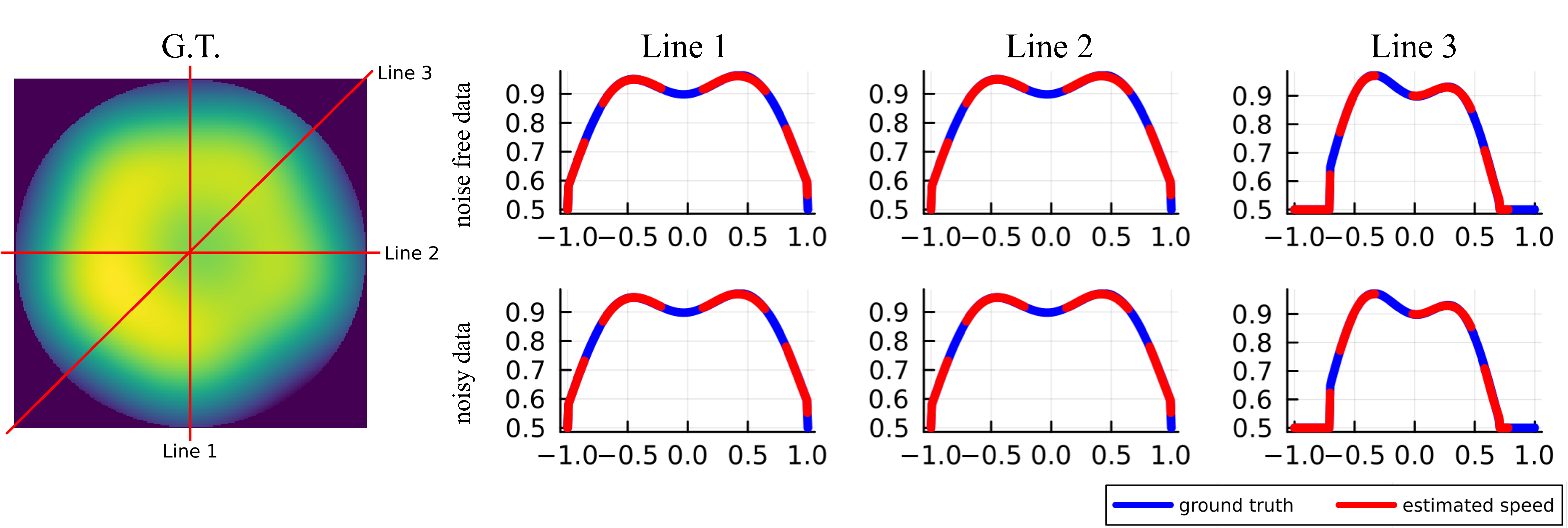}}
    \caption{In the left image, the ground truth is shown with three lines labeled Line 1, Line 2, and Line 3. The right images present cross-sectional views corresponding to these lines. The blue lines represent the ground truth, while the red dashed lines represent the estimated speed.}
    \label{fig:sliced_c_appendix}
\end{figure}

\section*{Acknowledgement}
This work was supported by the National Research Foundation of Korea grant funded by the Korea government(MSIT) (NRF-2022R1C1C1003464, RS-2023-00217116 and RS-2024-00333393).

\bibliographystyle{plain}

\end{document}